\documentclass[]{article}
\usepackage{amscd}
\usepackage{amsmath}
\usepackage{latexsym}
\usepackage{amsfonts}
\usepackage{amssymb}
\usepackage{amsthm}
\usepackage{bm}
\usepackage{graphicx,cases}
\usepackage{color,cite}

\newtheorem{lemma}{Lemma}[section]
\newtheorem{theorem}{Theorem}[section]

\newcommand{\beq}{\begin{equation}}
\newcommand{\eeq}{\end{equation}}
\newcommand{\ben}{\begin{eqnarray}}
\newcommand{\een}{\end{eqnarray}}
\newcommand{\beno}{\begin{eqnarray*}}
\newcommand{\eeno}{\end{eqnarray*}}

\theoremstyle{remark}
\newtheorem{remark}{Remark}[section]

\begin{document}
\title{Stability threshold for 2D shear flows of the\\ Boussinesq system near Couette}

\author{ \textsc{Dongfen Bian}  \\[1ex]
\normalsize School of Mathematics and Statistics\\
\normalsize Beijing Institute of Technology, Beijing 100081, China \\ 
\normalsize {biandongfen@bit.edu.cn}
\and
\textsc{Xueke Pu} \\[1ex] 
\normalsize School of mathematics and information science\\
\normalsize Guangzhou University, Guangzhou 510006, China\\ 
\normalsize {puxueke@gmail.com} \\
}
\date{}

\maketitle
\begin{abstract}
In this paper, we consider the  stability threshold for the shear flows of the Boussinesq system in a domain $\mathbb{T} \times \mathbb{R}$.  The main goal is to prove the nonlinear stability of the shear flow $(U^S,\Theta^S)=((e^{\nu t\partial_{yy}}U(y),0)^{\top},\alpha y)$ with $U(y)$ close to $y$ and $\alpha\geq0$. We separate two cases: one is $\alpha\geq 0$ small scaling with the viscosity coefficients and the case without smallness of $\alpha$ and fixed heat diffusion coefficient. The novelty here is that we don't require $\mu=\nu$ and only need to assume that $\mu$ is scaled with $\nu$ or fixed, where $\mu$ is the inverse of the Reynolds number and $\nu $ is the heat diffusion coefficient.
\end{abstract}
\begin{center}
 \begin{minipage}{120mm}
\end{minipage}
\end{center}
\begin{center}
 \begin{minipage}{120mm}
\end{minipage}
\end{center}

\section{Introduction}
\setcounter{section}{1}\setcounter{equation}{0}
In this paper, we consider the following 2D Boussinesq system
\begin{equation}\label{Eq1}
\begin{cases}
U_t+U\cdot \nabla U +\nabla P -\nu \Delta U=\Theta {e}_{2}, \\
\Theta_{t} + U\cdot \nabla\Theta-\mu\Delta\Theta = 0 , \\
\nabla \cdot {u} =0,
\end{cases}
\end{equation}
where $U= (U^{1},U^{2})^{\top}$ denotes the velocity, $\Theta$ the temperature, $P$ the pressure and $e_2=(0,1)^{\top}$. The spatial domain $\Omega$ here is taken to be $\Omega = \mathbb{T} \times \mathbb{R}$ with $\mathbb{T} = [0,2\pi]$ being the periodic box and $\mathbb{R}$ being the whole line, and $\mu$ is the inverse of the Reynolds number and $\nu$ is the heat conductivity.

Consider the following shear flow
\begin{equation}\label{shear}
U^S=(\bar{U}(t,y),0)^{\top}=(e^{\nu t\partial_{yy}}U(y),0)^{\top}, \ \ \ P^S=\alpha y^2/2,\ \ \ \Theta^S=\alpha y.
\end{equation}
It is easy to check that $(U^S,P^S,\Theta^S)$ is an exact solution to the Boussinesq equation \eqref{Eq1}. In this paper, we are concerned with shears $U(y)$ close to the Couette flow, i.e., shears that satisfy $\|\partial_yU-1\|_{H^{s}}\ll1$ for $s$ large enough. For the system \eqref{Eq1}, we consider the small perturbation of the shear flow profile
$$U(t,x,y)=U^S+\tilde u(t,x,y),\ \ P(t,x,y)=P^S+\tilde p(t,x,y), \ \ \Theta(t,x,y)=\Theta^S+\tilde\theta(t,x,y)$$
and the initial data for \eqref{Eq1} is set to be
$$U|_{t=0}=(U(y),0)^{\top}+\tilde u_{in}, \ \ \ \Theta|_{t=0}=\Theta^S+\tilde\theta_{in}.$$
When $\tilde u_{in}=\tilde\theta_{in}=0$, the solution reduces to the shear flow \eqref{shear}. The perturbation $(\tilde u,\tilde p,\tilde\theta)$ satisfy the following system
\begin{equation}\label{Per}
\begin{cases}
\partial_t{\tilde u}+\bar U\partial_x\tilde u+ \partial_y\bar U\left(
      \begin{array}{c}
       \tilde u^2 \\
        0 \\
      \end{array}
 \right)
+{\tilde u}\cdot\nabla{\tilde u} +\nabla\tilde p -\nu \Delta{\tilde u}=
\left(
\begin{array}{c}
0 \\
\tilde \theta \\
\end{array}
\right),\\
\nabla \cdot {\tilde u} =0,\\
\partial_t\tilde\theta + \bar U\partial_x\tilde \theta +\alpha \tilde u^2+ \tilde {u}\cdot \nabla\tilde \theta - \mu\Delta\tilde \theta = 0,
\end{cases}
\end{equation}
with initial data $\tilde u|_{t=0}=\tilde u_{in},\tilde \theta|_{t=0}=\tilde \theta_{in}$. We denote $\tilde \omega=\nabla^{\bot}\cdot \tilde u=-\partial_y\tilde u^1+\partial_x\tilde u^2$, then $\Delta\tilde \psi=\tilde \omega$, $\tilde u=\nabla^{\bot}\tilde \psi=
\left(
\begin{array}{c}
-\partial_y\Delta^{-1}\tilde \omega \\
\partial_x\Delta^{-1}\tilde \omega \\
\end{array}
\right).
$
Taking $\nabla^{\bot}\cdot$ to the system \eqref{Per}, we obtain
\begin{equation}\label{Per-W}
\begin{cases}
\partial_t{\tilde \omega}+\bar U\partial_x\tilde \omega -\partial_{yy}\bar U\partial_x\tilde \psi +\tilde {u}\cdot\nabla{\tilde \omega}-\nu\Delta{\tilde \omega} =\partial_x\tilde \theta,\\
\partial_t\tilde \theta + \bar U\partial_x\tilde \theta +\alpha\tilde u^2+ \tilde {u}\cdot \nabla\tilde \theta - \mu\Delta\tilde \theta = 0.
\end{cases}
\end{equation}
Consider a new coordinate system
\begin{equation}\label{Cortrans}
\left(
  \begin{array}{c}
    x \\
    y \\
  \end{array}
\right)
\mapsto
\left(
  \begin{array}{c}
    X \\
    Y \\
  \end{array}
\right)=
\left(
  \begin{array}{c}
    x-t\bar U(t,y) \\
    \bar U(t,y) \\
  \end{array}
\right),
\end{equation}
in which the new unknowns are denoted by $f(t,X,Y)=\tilde f(t,x,y)$. Then, the Boussinesq system \eqref{Per-W} reduces to the following system in the new coordinates
\begin{equation}\label{Per-new}
\begin{cases}
\partial_t\omega +{u}\cdot\nabla_t{\omega} =b\partial_X\psi +\nu\widetilde\Delta_t{\omega} +\partial_X\theta,\\
\partial_t\theta+ u\cdot\nabla_t\theta= \mu\widetilde\Delta_t\theta + (\mu-\nu)b\partial_Y^L\theta -\alpha u^2,\\
\omega=\Delta_t\psi,\ \ u=\nabla_t^{\bot}\psi,
\end{cases}
\end{equation}
where the operators $\nabla_t,\Delta_t,\widetilde\Delta_t$ and $\partial_Y^L$ will be given later on.

The main goal in this paper is to study the stability for the solution $(U^S,\Theta^S)=((\bar U(t,y),0)^{\top},\alpha y)$ for $\alpha\geq0$, i.e., the long time dynamics of the perturbation $(\tilde u,\tilde\theta)$ in the high Reynolds number regime $\nu\to0$. In particular, the main interest will focus on the stability threshold and determine how it scales with respect to $\nu$. Focus will also be put on the dependence of the heat diffusion coefficient $\mu\geq0$. The goal is, given a norm $\|\cdot\|_{X}$, to determine an index $\gamma=\gamma(X)$ such that
$$\|(\omega_{in},\theta_{in})\|_X\lesssim \nu^{\gamma}\text{ implies stability},$$
$$\|(\omega_{in},\theta_{in})\|_X\gg \nu^{\gamma}\text{ implies possible instability},$$
where $\mu$ is scaled with $\nu$ or fixed in this paper. Before we get into the details, we recall some recent mathematical results on the stability of shear flows.

When there is no heat convection, the Boussinesq system reduces to the Navier-Stokes equation, whose stability threshold was extensively studied very recently. For a given norm $\|\cdot\|_X$, the problem is then to determine an index $\gamma=\gamma(X)$ such that $\|\omega_{in}\|_X\lesssim \nu^{\gamma}$ implies stability and $\|\omega_{in}\|_X\gg \nu^{\gamma}$ implies possible instability. The $\gamma$ is sometimes referred to as the transition threshold. Nonlinear stability of the Couette flow in Sobolev space in the case of infinite channel was considered first by Romanov \cite{Rom73}. For the 2D Couette flow of the Navier-Stokes equation, Bedrossian \emph{et al} \cite{BMV16,BGM19} proved the Couette flow is uniformly stable at high Reynolds number and there is no subcritical transition in Gevrey-$m$ with $m<2$, when the initial perturbation is sufficiently small. It was also shown that the Couette flow is nonlinearly stable for the 2D Euler equations for sufficiently smooth Gevrey perturbations \cite{BM13}. In the case of 3D Couette flow for the Navier-Stokes equation, it was shown that $\gamma=1$ in Gevrey-$m$ for $m<2$ by Bedrossian \emph{et al} \cite{BGM15a,BGM15b} and $\gamma\leq 3/2$ in $X=H^s$ for $s>7/2$ by Bedrossian \emph{et al} \cite{BGM17}, which is consistent with the numerical estimation of $31/20$ given in Reddy et al \cite{Reddy}. When the shear flow is only `near' the Couette flow in some sense, Bedrossian \emph{et al} \cite{BVW18} showed the Sobolev stability threshold is $\gamma=1/2$, i.e., if the data is initially $\varepsilon$-close to the shear flow, then the solution remains $\varepsilon$-close to the shear flow. Zillinger \cite{Zill17} and Jia \emph{et al} \cite{Jia20,Jia20SIAM} also studied the linear inviscid damping for the monotone shear flow. See their papers for exact statement. Wei \emph{et al} \cite{WZZ18,WZZ20,WZ18} obtained the optimal decay estimate for the linearized problem around monotone shear flows and studied the linear inviscid damping in the case of nonmonotone shear flows. We would like to remark that Grenier \emph{et al}  \cite{G-G-T16} proved the spectral instability of general symmetric shear flows of the Navier-Stokes equations at a high Reynolds number in a 2D channel.

For the Boussinesq system, the study of stability for shear flows has a long history, starting from the work of Taylor \cite{Tay}, Goldstein \cite{Gold} and Synge \cite{Synge} for the linearized 2D Boussinesq system. It also arouses interests for many authors in recent years. For local and global well-posedness with/without viscosity, one may refer to \cite{BZD2005,CD81,Chae06,CN97,HL05} and the references therein. The stability and enhanced dissipation phenomenon for the linearized 2D Boussinesq system was recently studied by Lin and Zeng \cite{LZ11} for the case $U(y)=y$ and $\alpha=0$ for stratified fluids, and by Tao and Wu \cite{TW19} by using the hypocoercivity method introduced in Villani's paper \cite{Vil09}. Nonlinear stability and large time behaviors of Couette flow $(U^S,\Theta^S)=((y,0)^{\top},0)$ for the 2D Boussinesq system were obtained by Deng \emph{et al} \cite{DWZ2004}, for both the full dissipation case and the case with only vertical dissipation, when $\mu=\nu$. In particular, a new time independent bounded Fourier multiplier different from that in \cite{BGM17,BVW18} was introduced, enabling the authors to extract the enhanced dissipation. In the stratified case, the linear inviscid damping is studied by Yang and Lin \cite{YL18} for the Couette flow and for shear flows near Couette in the 2D stably stratified regime by Bianchini \emph{et al} \cite{BZD2005}. 
Zillinger \cite{Zill20} studied the stability for the linearized system with partial dissipation and obtained nonlinear stability for the full dissipated 2D Boussinesq system for the shear flow $(U^S,\Theta^S)=((\beta y,0)^{\top},\alpha y)$ for $\alpha\geq0$. See also Zillinger \cite{Zill2011} for the case that $\Theta^S=\alpha y$ is replaced with a general $\Theta^S=T(y)$ with suitable restrictions, and Doering \emph{et al} \cite{DWZZ18} for long time behavior of the 2D Boussinesq equations without buoyancy diffusion. When there is no heat diffusion $\mu=0$, Masmoudi \emph{et al} \cite{MSZ20} studied the stability of Couette flow for the Navier-Stokes Boussinesq system for the initial perturbation in Gevrey-$\frac1s$ for $1/3<s\leq 1$ in the domain $\Bbb T\times\Bbb R$.

In this paper, we will focus ourselves on the stability of shear flows that is close to the Couette flow with hydrostatic balance. That is to say, we will mainly focus on the shear flow $(U^S,\Theta^S)=((e^{\nu t\partial_{yy}}U(y),0)^{\top},\alpha y)$ with $U(y)$ close to $y$ and $\alpha\geq0$.  The main results are stated in the following two theorems.
\begin{theorem}[$\alpha\ll1$]\label{thm1}
Let $N\geq5$, $\nu,\mu\in(0,1)$ with $\nu\lesssim \mu$, and
\begin{equation*}
\|U'-1\|_{H^s}+\|U''\|_{H^s}\leq \delta\ll 1
\end{equation*}
for some small $\delta$ independent of $\nu$ and $\mu$. There exists $\gamma_1$ and $\gamma_2$ such that if $$\varepsilon_1=\|\omega_{in}\|_{H^N}\leq \gamma_1\min(\nu,\mu)^{1/2},\ \ \ \varepsilon_2=\|\theta_{in}\|_{H^N}\leq \gamma_2\sqrt{\nu\mu}\min(\nu,\mu)^{1/2},$$
and $0\leq \alpha<\nu^{1/3}\mu^{1/2}\varepsilon_2/ \varepsilon_1$, then the solution with initial data $(\omega_{in},\theta_{in})$ is global in time and for any $T>0$, it holds
\begin{equation}\label{eq8}
\begin{split}
\|\omega\|_{L^\infty(0,T; H^N)}+\nu^{1/2}\|\nabla_L\omega\|_{L^2(0,T;H^N)}\lesssim & \varepsilon_1,\\
\|\theta\|_{L^\infty(0,T; H^N)}+\mu^{1/2}\|\nabla_L\theta\|_{L^2(0,T;H^N)}\lesssim & \varepsilon_2,\\
\|u^X_0\|_{L^{\infty}(0,T;L^2)} +\nu^{1/2}\|\partial_Yu^X_0\|_{L^2(0,T;L^2)} \lesssim & \varepsilon_1,
\end{split}
\end{equation}
and
\begin{equation}\label{eq9}
\begin{split}
\|\omega_{\neq}\|_{L^2(0,T;H^N)} \lesssim\varepsilon_1\nu^{-1/6},\ \ \ \|\theta_{\neq}\|_{L^2(0,T;H^N)} \lesssim\varepsilon_2\nu^{-1/6},
\end{split}
\end{equation}
where the implicit constants do not depend on $\nu,\mu$ or the initial data.
\end{theorem}
\begin{remark}
In the above theorem, unlike previous papers of Zillinger \cite{Zill20} or Deng \emph{et al } \cite{DWZ2004}, we do not require $\nu=\mu$. This means that we can let both $\nu$ and $\mu$ tend to zero independently. In particular, if $\nu\sim\mu$, this theorem implies that when the shear flow is near the Couette flow in some sense, then the solution is nonlinearly stable provided $\|\omega_{in}\|_{H^{N}}\lesssim \nu^{1/2}$, in accordance with the index $\gamma=1/2$ in \cite{BMV16} in the Navier-Stokes equation case.
\end{remark}

The above theorem excludes the case when $\alpha$ is not small. In this case, we have the following theorem.
\begin{theorem}[$\alpha$ not small]\label{thm2}
Let $N\geq5$, $\nu\in(0,1)$, $\mu\geq2$, $\alpha>0$ be fixed and
\begin{equation*}
\|U'-1\|_{H^s}+\|U''\|_{H^s}\leq \delta\ll 1
\end{equation*}
for some small $\delta$ independent of $\nu$ and $\mu$. Assume
$$\alpha\|\omega_{in}\|^2_{H^N}= \varepsilon_1^2\leq \varepsilon^2\lesssim \gamma_1\nu,\ \ \ \|\nabla\theta_{in}\|_{H^N}^2\leq \varepsilon^2\lesssim \gamma_1\nu,$$
then the solution with initial data $(\omega_{in},\theta_{in})$ is global in time and satisfies the following estimates. For any $T>0$,
\begin{equation}\label{eq27}
\begin{split}
\alpha\|A\omega(T)\|^2_{L^2}  & +\|\nabla_LA\theta\|^2_{L^2}+ \int_0^T\left(\nu\alpha\|\nabla_LA\omega\|_{L^2}^2 +\frac{\alpha}{2}\|\sqrt{-\dot MM}\omega\|^2_{H^N}\right)dt\\
& +\int_0^T\left(\|\sqrt{-\Delta_L}\sqrt{-\dot MM}\theta\|_{H^N}  +\frac{\mu}{4}\|\Delta_LA\theta\|_{L^2}^2\right)dt \lesssim \varepsilon^2
\end{split}
\end{equation}
and
\begin{equation}\label{eq26}
\begin{split}
\|\omega_{\neq}\|_{L^2(0,T;H^N)} \lesssim\varepsilon\alpha^{-1/2}\nu^{-1/6},\ \ \ 
\end{split}
\end{equation}
where the implicit constants do not depend on $\nu,\mu$ or the initial data.
\end{theorem}

\begin{remark}
Unlike the small $\alpha$ case, we cannot estimate the two unknowns $\omega$ and $\theta$ separately, rather a combination of them, since in this large $\alpha$ case, we need to cancel the gravitational term with the $\alpha$ term in the heat equation at different energy level, with proper sign $\alpha>0$. However, a new linear term $\langle \partial_X\partial_Y^LA\theta ,A\theta\rangle$ coming from the time derivative term will appear and cannot be bounded generally without assuming $\mu$ large suitably.
\end{remark}
\begin{remark}
As in the small $\alpha$ case, again we do not require $\nu=\mu$. In the proof process, we need to assume $\nu\leq 2\mu$ to control the linear term $b\partial_Y^L\theta$ in the heat equation, which holds automatically from the assumption $\mu\geq2$.
\end{remark}


The above two theorems will be proved in the subsequent section. First, we will make the change of coordinates and hence the dependent variables. Then we will prove the two theorems by continuity method. Before we get into the next section, we list some frequently used notations. We write $\langle x\rangle=\sqrt{1+x^2}$. For a function $f(x,y)$, denote by
$$f_0(y)=\int_{\Bbb T}f(x,y)dx,\ \ \ \ f_{\neq}(x,y)=f(x,y)-f_0(y),$$
the projection onto the zero and nonzero frequencies, respectively, w.r.t. $x$ variable.



\section{Nonlinear Stability}
\setcounter{section}{2}\setcounter{equation}{0}

\subsection{Change of dependent variables}
Consider the more general class of shear flows with initial data $(U(y),0)$ which is sufficiently close to Couette flow, in the sense that
\begin{equation}\label{eq24}
\|U'-1\|_{H^s}+\|U''\|_{H^s}\leq \delta\ll 1
\end{equation}
for some small $\delta$ independent of $\nu$ and $\mu$, and for some $s\geq2+N$. From the definition of $\bar U(t,y)=e^{\nu t\partial_{yy}}U(y)$, which satisfies a 1D heat equation $\partial_tu=\nu\partial_{yy}u$ with initial data $u(t=0,\cdot)=U(\cdot)$, one immediately has
\begin{equation}\label{eq25}
\begin{split}
\sup_{t>0}\|\bar U'(t,\cdot)-1\|_{H^s}\leq & \|U'-1\|_{H^s},\\
\sup_{t>0}\|\bar U''(t,\cdot)\|_{H^s}\leq & \|U''\|_{H^s},\\
\|\bar U''(\cdot,\cdot)\|_{L^2_tH^s_y}\leq & \delta\nu^{-1/2}.
\end{split}
\end{equation}

Let us consider the following coordinate transform
\begin{equation}\label{Cortrans}
\left(
  \begin{array}{c}
    x \\
    y \\
  \end{array}
\right)
\mapsto
\left(
  \begin{array}{c}
    X \\
    Y \\
  \end{array}
\right)=
\left(
  \begin{array}{c}
    x-t\bar U(t,y) \\
    \bar U(t,y) \\
  \end{array}
\right).
\end{equation}
Denote $a(t,Y)=\partial_y\bar U(t,y)$ and $b(t,Y)=\partial_{yy}\bar U(t,y)$. One can deduce that
$\partial_x=\partial_X$ and $\partial_y=\partial_y\bar U(t,y)(\partial_Y-t\partial_X)=a(t,Y)(\partial_Y-t\partial_X)$. Obviously, one has $$b(t,Y)=\partial_y(a(t,Y))=\partial_Ya\cdot \partial_yY=\partial_Ya\cdot\partial_y\bar U(t,y)=a\partial_Ya.$$
From the above smallness condition, for $\delta$ sufficiently small, there holds
\begin{equation*}
\|\partial_Y(a^2-1)\|_{L^2H^{\sigma}_y}=\|2b\|_{L^2H^{\sigma}_y}\lesssim\|\bar U''\|_{L^2H^{\sigma}_y} \lesssim \delta\nu^{-1/2},
\end{equation*}
and similarly
\begin{equation}\label{eq3}
\|a-1\|_{H^{\sigma}} +\|b\|_{H^{\sigma}} \leq 2\delta,
\end{equation}
for $\sigma>2$ and $\delta>0$ small enough.

To distinguish the old and new coordinates, we get rid of the $\tilde\cdot$ above a function and set $h(t,X,Y)=\tilde h(t,x,y)$, where the shorthands $X=X(t,x,y)$ and $Y=Y(t,x,y)$ are used. In the new coordinates, differential operators are modified as
$$\nabla\tilde h(t,x,y)=
\left(
\begin{array}{c}
\partial_x\tilde h \\
\partial_y\tilde h \\
\end{array}
\right)
=\left(
   \begin{array}{c}
     \partial_Xh \\
     a(\partial_Y-t\partial_X)h \\
   \end{array}
 \right)
=\left(
   \begin{array}{c}
     \partial^t_Xh \\
     \partial^t_Yh \\
   \end{array}
 \right)
=:\nabla_th,
$$
and the operator $\Delta\mapsto \Delta_t$ is defined via
$$\Delta \tilde h(t,x,y)=(\partial_X^2+a^2\partial^{L}_{YY} +b\partial^L_Y)h=:\Delta_th.$$
In particular, in the Couette flow case $U(y)=y$, we have $\bar U(t,y)=y$. In this case, the transform of coordinate \eqref{Cortrans} and the operators in the new coordinates both can be simply computed. Indeed we have $a\equiv1$ and $b\equiv0$ and the operator $\nabla_t$ reduces to
$$
\nabla_L=\left(
   \begin{array}{c}
     \partial_X \\
     \partial^L_Y \\
   \end{array}\right)=\left(
   \begin{array}{c}
     \partial_X \\
     \partial_Y-t\partial_X \\
   \end{array} \right)
$$
and $\Delta_t$ reduces to $\Delta_L=\partial_{X}^2+\partial_{YY}^L= \partial_{X}^2+(\partial_{Y}-t\partial_X)^2$. In particular, when $\bar U(t,y)=y$, $\partial_{yy}\bar U=0$ and the term $\partial_{yy}\bar U\partial_x\tilde\psi$ in \eqref{Per-W} vanishes. However, In the shear flow near Couette case, $\bar U(t,y)=e^{\nu t\partial_{yy}}U(y)$, the transform depends on time and some functions involved with $\bar U(t,y)$ need be introduced.

Thanks to $\partial_t\bar U=\nu\partial_{yy}\bar U$, the time derivative of $\tilde h$ reads
$$\partial_t\tilde h=\partial_th-Y\partial_Xh+\nu b(t,Y)(\partial_Y-t\partial_X)h =\partial_th-Y\partial_Xh+\nu b\partial_Y^Lh.$$
The relationship between $\omega$ and $\psi$ is given through
\begin{equation*}
\omega=\Delta_t\psi =(\partial_X^2+a^2\partial^{L}_{YY}+b\partial^L_Y)\psi =\Delta_L\psi +\left((a^2-1)\partial_{YY}^L+b\partial^L_Y\right)\psi.
\end{equation*}
If we further modify the Laplacian operator as
$$\widetilde\Delta_t=\Delta_L +(a^2-1)\partial_{YY}^L,$$
then
$$\omega=\widetilde\Delta_t\psi+b\partial^L_Y\psi.$$

Define
\begin{equation*}
\begin{split}
\omega(t,X,Y)=&\tilde\omega(t,x,y),\ \ \
\theta(t,X,Y)=\tilde\theta(t,x,y),\\
\psi(t,X,Y)=&\tilde\psi(t,x,y),\ \ \
u(t,X,Y)=\tilde u(t,x,y).
\end{split}
\end{equation*}
Then, the Boussinesq system \eqref{Per-W} reduces to the following system in the new coordinates
\begin{equation}\label{Per-new}
\begin{cases}
\partial_t\omega +{u}\cdot\nabla_t{\omega} =b\partial_X\psi +\nu\widetilde\Delta_t{\omega} +\partial_X\theta,\\
\partial_t\theta+ u\cdot\nabla_t\theta= \mu\widetilde\Delta_t\theta + (\mu-\nu)b\partial_Y^L\theta -\alpha u^2,\\
\omega=\Delta_t\psi,\ \ u=\nabla_t^{\bot}\psi.
\end{cases}
\end{equation}
We remark that when $\mu=\nu$, the linear term $(\mu-\nu)b\partial_Y^L\theta$ vanish and causes no trouble. When $\mu\neq\nu$, this term is linear and generally causes growth. But since $b$ is small, in the case $0\leq\alpha\ll1$, this term is very small and harmless if we require $|\mu-\nu|\mu^{-1}\lesssim 1$, which can be guaranteed by letting $\nu\lesssim\mu$. In the case that $\alpha$ is not small, this term can be absorbed by the diffusion term when $\nu\leq 2\mu$, which holds naturally as we have assumed $\mu\geq2$ in this case. For the linear term $\alpha u^2$, it will be canceled with the linear gravitational term $\partial_X\theta$ when we take energy estimates at different Sobolev level with proper sign $\alpha>0$.

We also remark that since $\bar U(t,y)$ is independent of $x$, from the definition of the zero and nonzero frequencies, it holds that
$$(\nabla_th)_0=\nabla_th_0,\ \ \ \ (\nabla_th)_{\neq}=\nabla_th_{\neq}.$$
And by the divergence free condition we deduce that $u_0^X=-a\partial_Y\Delta_t^{-1}\omega_0$ and the $X$-average of the second velocity component vanishes, i.e., $u_0^Y=0$, thanks to the fact that $a(t,Y)$ and $b(t,Y)$ are independent of $X$.

From the assumptions \eqref{eq24} and \eqref{eq25}, we get the following Lemma.
\begin{lemma}
Assume \eqref{eq24} and \eqref{eq25} hold, then it holds for $\sigma>2$ and $\delta>0$ sufficiently small that
\begin{equation}\label{eq5}
\begin{split}
\|\partial_Y(a^2-1)\|_{L^2H_y^{\sigma}}\lesssim \|\bar U''\|_{L^2H_y^{\sigma}}\lesssim \delta \nu^{-1/2}
\end{split}
\end{equation}
and
\begin{equation}\label{eq6}
\begin{split}
\|a-1\|_{H^{\sigma}} +\|b\|_{H^{\sigma}} \leq 2\|\bar U'-1\|_{H^{\sigma}}+2\|\bar U''\|_{H^{\sigma}} \lesssim 2\delta.
\end{split}
\end{equation}
\end{lemma}

Before we state our main results, we would like to list some useful inequalities involved with  the multiplier $M(t,k,\xi)$ introduced in \cite{BGM17,BVW18} and some helpful operators.
\begin{lemma}
There exists a Fourier multiplier $M(t,k,\xi)$ with the following properties:
\begin{equation*}
\begin{split}
M(0,k,\xi)=&M(t,0,\xi)=1,\\
1\geq M(t,k,\xi)\geq & c,\\
\frac{\dot M}{M}\geq & \frac{|k|}{k^2+|\xi-kt|^2}\ \ \ \text{for} \ k\neq 0,\\
\frac{\partial_xM(k,\xi)}{M(k,\xi)}\lesssim & \frac{1}{|k|}\ \ \text{for}\ k\neq0,\ \text{uniformly in} \ \xi,\\
1\lesssim & \nu^{-1/6}\left(\sqrt{-\dot MM(t,k,\xi)}+\nu^{1/2}|k,\xi-kt|\right)\ \ \text{for }\ k\neq 0,\\
\sqrt{-\dot MM(t,k,\xi)}\lesssim & \langle\eta-\xi\rangle\sqrt{-\dot MM(t,k,\xi)},
\end{split}
\end{equation*}
for some constant $c\in (0,1)$ which is independent of $\nu$.
\end{lemma}
One may refer to \cite{BGM17,BMV16} for detailed proofs. From these inequalities, it is easy to deduce that
\begin{equation}\label{eq15}
\|f_{\neq}\|_{L^2H^N}\lesssim \nu^{-1/6}\left(\|\sqrt{-\dot MM}f_{\neq}\|_{L^2H^N}+\nu^{1/2}\|\nabla_Lf_{\neq}\|_{L^2H^N}\right)
\end{equation}
holds for any $f$ and $N\geq0$. Let $f\in H^N$ and $N>1$, then for $\delta$ sufficiently small, there holds
$$\|\sqrt{-\dot MM}\Delta_L\Delta_t^{-1}f_{\neq }\|_{H^N}\lesssim \|\sqrt{-\dot MM}f_{\neq }\|_{H^N}.$$
This implies that $\Delta_L\Delta_t^{-1}$ is almost identity and indeed for $\delta$ sufficiently small, it holds that
$$\|\Delta_L\Delta_t^{-1}f_{\neq }\|_{H^N}\lesssim \|f_{\neq }\|_{H^N}.$$
Moreover, using the above inequality and Parseval theorem, when $\delta$ is sufficiently small, we can get for $f\in H^N$ with $N>1$
\begin{equation*}
\begin{split}
\|\partial_X\Delta_t^{-1}f_{\neq} \|_{H^N}= & \|(\partial_X\Delta_L^{-1}\Delta_L \Delta_t^{-1}f_{\neq} \|_{H^N}\leq \|\partial_X\Delta_L^{-1}f_{\neq} \|_{H^N}\\
= & \|k(k^2+|\xi-kt|^2)^{-1}\widehat{f_{\neq}} \|_{H^N}\leq \|\widehat{f_{\neq}} \|_{H^N} =\|{f_{\neq}} \|_{H^N}.
\end{split}
\end{equation*}
In other words, in the operator norm, we have $\|\partial_X\Delta_t^{-1}\|\lesssim 1$ in $H^N$ for $N>1$.

\subsection{Proof of Theorem \ref{thm1}}

In what follow, we set $A=M\langle D\rangle^N$, where $\langle D\rangle=\sqrt{1+D^2}$ and $D=i\nabla$ is the Fourier multiplier $(k,\xi)^{\bot}$. From the Plancherel's inequality, it holds that
$$\|f\|_{H^N}\lesssim \|A(t)f\|_{L^2}\lesssim \|f\|_{H^N},$$
and $\|A(0)f\|_{L^2}=\|f\|_{H^N}$.

We will employ the continuity method to prove the theorem, based on establishing the following \emph{a priori} estimates
\begin{equation}\label{assumption}
E_{\omega}(T)\leq 8\varepsilon_1^2,\ \ \ \ E_{\theta}(T)\leq 8\varepsilon_2^2,
\end{equation}
and
\begin{equation}\label{eq2}
\|u^X_0\|_{L^{\infty}(0,T;L^2)} +\nu^{1/2}\|\partial_Yu^X_0\|_{L^2(0,T;L^2)} \leq 8\varepsilon_1,
\end{equation}
where
\begin{equation*}
E_{\omega}(T):=\|A\omega\|^2_{L^{\infty}(0,T;L^2)} +\nu\|\nabla_LA\omega\|^2_{L^2(0,T;L^2)} +\|\sqrt{-\dot MM}\langle D\rangle^N\omega\|^2_{L^2(0,T;L^2)},
\end{equation*}
and
\begin{equation*}
E_{\theta}(T):=\|A\theta\|^2_{L^{\infty}(0,T;L^2)} +\mu\|\nabla_LA\theta\|^2_{L^2(0,T;L^2)} +\|\sqrt{-\dot MM}\langle D\rangle^N\theta\|^2_{L^2(0,T;L^2)}.
\end{equation*}

In particular, the first two inequalities of \eqref{eq8} are implied by \eqref{assumption}, and \eqref{eq9} is implied by \eqref{eq15} and \eqref{assumption}. From local well-posedness, there exists some small $T>0$ such that \eqref{assumption} and \eqref{eq2} hold with `8' replaced by `2' and all the quantities in $E_{\omega}$ and $E_{\theta}$ are continuous in time. To employ the continuity method, we define $T_*\leq\infty$ to be the maximal time such that
\begin{equation}\label{assump}
E_{\omega}(T_*)\leq 8\varepsilon_1^2,\ \ \ \ E_{\theta}(T_*)\leq 8\varepsilon_2^2,
\end{equation}
In particular $T_*>T$. This is the bootstrap assumption. Next we will prove that the inequalities in \eqref{assump} and \eqref{eq2} indeed hold with `8' replaced with `4' on the RHS. Therefore, the maximal time has to be $\infty$, concluding the proof, by continuity method.


Now, let $T>0$ be such that \eqref{assumption} holds. Applying the operator $A$ to \eqref{Per-new}, and then taking inner product of $A\omega$ and $A\theta$ respectively to obtain the following
\begin{equation}\label{eq20}
\begin{split}
\frac12\|A&\omega(T)\|^2_{L^2} +\nu\|\nabla_LA\omega\|^2_{L^2_TL^2} +\|\sqrt{-\dot MM}\langle D\rangle^N\omega\|^2_{L^2_TL^2} \\
=& \frac12\|A\omega_{in}\|^2_{L^2} -\int_0^T\int_{\Bbb R}\int_{\Bbb T} A(u\cdot\nabla_t\omega)A\omega dXdYdt +\int_0^T\int_{\Bbb R}\int_{\Bbb T} A(b\partial_X\psi)A\omega dXdYdt \\
&+\nu\int_0^T\int_{\Bbb R}\int_{\Bbb T}  A((a^2-1)\partial^L_{YY}\omega)A\omega dXdYdt + \int_0^T\int_{\Bbb R}\int_{\Bbb T} A\partial_X\theta A\omega dXdYdt\\
=& : \frac12\|A\omega_{in}\|^2_{L^2} -\mathcal T_{\omega} +\cal S+\cal D_{\omega}+\cal{T}_{\omega\theta},
\end{split}
\end{equation}
and
\begin{equation}\label{eq21}
\begin{split}
\frac12\|A&\theta(T)\|^2_{L^2} +\mu\|\nabla_LA\theta\|^2_{L^2_TL^2} +\|\sqrt{-\dot MM}\langle D\rangle^N\theta\|^2_{L^2_TL^2} - \frac12\|A\theta_{in}\|^2_{L^2}\\
=& -\int_0^T\int_{\Bbb R}\int_{\Bbb T}  A(u\cdot\nabla_t\theta)A\theta dXdYdt +\mu\int_0^T\int_{\Bbb R}\int_{\Bbb T}  A((a^2-1)\partial^L_{YY}\theta)A\theta dXdYdt \\
&+(\mu-\nu)\int_0^T\int_{\Bbb R}\int_{\Bbb T}  A(b\partial_Y^L\theta)A\theta dXdYdt -\alpha\int_0^T\int_{\Bbb R}\int_{\Bbb T}  A(\partial_X\Delta_t^{-1}\omega)A\theta dXdYdt\\
=& : \frac12\|A\theta_{in}\|^2_{L^2} -{\cal T}_{\theta}  +{\cal D}_{\theta} +{\cal T}_{b} -{\cal T}_{\theta\omega}.
\end{split}
\end{equation}
It can be shown that
\begin{equation*}
{\cal T}_{\omega}\lesssim \varepsilon_1^3\nu^{-1/2},\ \ \ {\cal D}_{\omega}\leq \delta\varepsilon_1^2,\ \ \ {\cal T}_{\theta}\lesssim \varepsilon_1\varepsilon_2^2\nu^{-1/2},\ \ \ {\cal D}_{\theta}\leq \delta\varepsilon_2^2
\end{equation*}
and the term $\cal S$ can be absorbed on the left hand side.

Some of the estimates of ${\cal T}_{\omega}$ and ${\cal D}_{\omega}$ are similar to that in \cite{BVW18} and will only be sketched in the following. First, we consider ${\cal T}_{\theta}$ and split the velocity field into two parts according to the zero and the nonzero modes, i.e.,
\begin{equation*}
\begin{split}
{\cal T}_{\theta}=\int_0^T\int_{\Bbb R}\int_{\Bbb T}  A(u^X_0\partial_X\theta)A\theta dXdYdt + \int_0^T\int_{\Bbb R}\int_{\Bbb T}  A(\nabla^{\bot}_t\Delta_t^{-1}\omega\cdot\nabla_t\theta)A\theta dXdYdt=:{\cal T}_{\theta,0}+{\cal T}_{\theta,\neq}.
\end{split}
\end{equation*}
For the zero part, we have
\begin{equation*}
\begin{split}
{\cal T}_{\theta,0}= & \int_0^T\int_{\Bbb R}\int_{\Bbb T}  A(u^X_0\partial_X (\theta_0+\theta_{\neq}))A(\theta_0+\theta_{\neq}) dXdYdt\\
= &\int_0^T\int_{\Bbb R}\int_{\Bbb T}  A(u^X_0\partial_X \theta_{\neq})A \theta_{\neq} dXdYdt\\
= &\int_0^T\int_{\Bbb R}\int_{\Bbb T}  \left(A(u^X_0\partial_X \theta_{\neq})-u_0^X\partial_XA\theta_{\neq}\right)A \theta_{\neq} dXdYdt.
\end{split}
\end{equation*}
Here we have used in the second line the fact that
\begin{equation*}
\begin{split}
\int_0^T\int_{\Bbb R}\int_{\Bbb T}  A(u^X_0\partial_X \theta_{\neq})A \theta_0 dXdYdt = -\int_0^T\int_{\Bbb R}\int_{\Bbb T}  A\left([a\partial_Y^L\Delta_t^{-1}]\partial_X \theta_{\neq}\right)A \theta_0 dXdYdt=0,
\end{split}
\end{equation*}
by integration by parts and the fact that $a$ and $b$ are both independent of $X$, and in the third line
\begin{equation*}
\begin{split}
\int_0^T\int_{\Bbb R}\int_{\Bbb T}  u_0^X\partial_XA\theta_{\neq} A \theta_{\neq} dXdYdt=0,
\end{split}
\end{equation*}
both by integration by parts and the divergence free condition. By Plancherel's theorem,
\begin{equation}\label{eq4}
\begin{split}
{\cal T}_{\theta,0}\sim \sum_{k\neq0}\int_0^T\int_{\Bbb R}\int_{\Bbb R} \left(A(k,\xi)-A(k,\xi-\eta)\right)\widehat{u_0^X}(\eta)ik \widehat{\theta}(k,\xi-\eta)d\eta A(k,\xi)\overline{\widehat\theta}(k,\xi) d\xi dt
\end{split}
\end{equation}
by a universal constant. A direct computation gives
\begin{equation*}
\begin{split}
\left|k(A(k,\xi)-A(k,\xi-\eta))\right|\lesssim \left((1+k^2+(\xi-\eta)^2)^{N/2} +(1+k^2+\xi^2)^{N/2} \right)|\eta|.
\end{split}
\end{equation*}
Recall $u_0^X=-a\partial_Y\Delta_t^{-1}\omega_0$ and hence
$$\|\partial_Y\left(u_0^X\right)\|_{H^N}\lesssim \|\partial_Y(a-1)\|_{H^N} \|\partial_Y\Delta_t^{-1}\omega_0\|_{H^N}+ \|a\|_{H^N} \|\partial_{YY}\Delta_t^{-1}\omega_0\|_{H^N}.$$
By Young's convolution inequality and the inequality
$$\|\widehat f\|_{L^1}\lesssim \|f\|_{H^N},\ \ \ N>1,$$
we have from \eqref{eq4} that
\begin{equation}\label{eq7}
\begin{split}
{\cal T}_{\theta,0}\lesssim & \|\partial_Y(a\partial_Y\Delta_t^{-1}\omega_0) \|_{L^{\infty}H^N} \|\theta_{\neq }\|_{L^2H^N} \|A\theta_{\neq }\|_{L^2L^2}\\
\lesssim & \left(\|\partial_Y(a-1)\|_{L^{\infty}H^N} \|\partial_Y\Delta_t^{-1}\omega_0\|_{L^{\infty}H^N}+ \|a\partial_{YY}\Delta_t^{-1}\omega_0\|_{L^{\infty}H^N}\right)\\
& \cdot \|\theta_{\neq }\|_{L^2H^N} \|A\theta_{\neq }\|_{L^2L^2}\\
\lesssim & \left(\delta\|\partial_Y\Delta_t^{-1} \omega_0\|_{L^{\infty}H^N} +\|\partial_{YY}\Delta_t^{-1} \omega_0\|_{L^{\infty}H^N}\right)\|\theta_{\neq }\|_{L^2H^N} \|A\theta_{\neq }\|_{L^2L^2}\\
\lesssim & \left(\|u_0^X\|_{L^{\infty}L^2} +\|\omega_0\|_{L^{\infty}H^N}\right)\|\theta_{\neq }\|_{L^2H^N} \|A\theta_{\neq }\|_{L^2L^2},
\end{split}
\end{equation}
where we have used \eqref{eq6},
\begin{equation*}
\begin{split}
\|a\partial_{YY}\Delta_t^{-1} \omega_0\|_{H^N}
\lesssim & \|a\|_{L^{\infty}}\|\partial_{YY}\Delta_t^{-1} \omega_0\|_{H^N} +\|\nabla a\|_{H^{N-1}}\|\partial_{YY}\Delta_t^{-1} \omega_0\|_{H^N}\\
\lesssim & (\|a-1\|_{H^N}+1)\|\partial_{YY}\Delta_t^{-1} \omega_0\|_{H^N}\\
\lesssim & \|\partial_{YY}\Delta_t^{-1} \omega_0\|_{H^N},
\end{split}
\end{equation*}
and
\begin{equation*}
\begin{split}
\|\partial_{Y}\Delta_t^{-1} \omega_0\|_{H^N} \lesssim & \|\partial_{Y}\Delta_t^{-1} \omega_0\|_{L^2} +\|\partial_{YY}\Delta_t^{-1} \omega_0\|_{H^N}\\
\lesssim & \|u_0^X\|_{L^2} +\|\omega_0\|_{H^N}.
\end{split}
\end{equation*}
From \eqref{eq7}, \eqref{eq2} and \eqref{eq9}, one immediately has
\begin{equation}\label{eq10}
\begin{split}
|{\cal T}_{\theta,0}|\lesssim \varepsilon_1\varepsilon_2^2\mu^{-1/3}.
\end{split}
\end{equation}

For the nonzero mode ${\cal T}_{\theta,\neq}$, we have by the smallness condition \eqref{eq3}
\begin{equation}\label{eq11}
\begin{split}
|{\cal T}_{\theta,\neq}|\lesssim  & \|\nabla_t^{\bot}\Delta_t^{-1}\omega_{\neq}\|_{L^2H^N} \|\nabla_t\theta_{\neq}\|_{L^2H^N}  \|A\theta\|_{L^\infty L^2}\\
\lesssim & \|\nabla_L\Delta_t^{-1}\omega_{\neq}\|_{L^2H^N} \|\nabla_L\theta_{\neq}\|_{L^2H^N}  \|A\theta\|_{L^\infty L^2}\\
\lesssim & \|\sqrt{-\dot MM}\omega_{\neq}\|_{L^2H^N} \|\nabla_L\theta_{\neq}\|_{L^2H^N}  \|A\theta\|_{L^\infty L^2}\lesssim \varepsilon_1\cdot\varepsilon_2\mu^{-1/2}\cdot \varepsilon_2
\end{split}
\end{equation}
thanks to the fact that $H^N$ is an algebra for $N>1$ and Lemma 4.3 in \cite{BVW18}. Combining \eqref{eq10} and \eqref{eq11}, one then has (recall $\nu,\mu<1$)
\begin{equation}\label{eq16}
\begin{split}
|{\cal T}_{\theta}|\lesssim  \varepsilon_1\varepsilon_2^2\mu^{-1/2}.
\end{split}
\end{equation}
Similar treatments leads to the estimate for ${\cal T}_{\omega}$,
\begin{equation}\label{eq17}
\begin{split}
|{\cal T}_{\omega}|\lesssim  \varepsilon_1^3\nu^{-1/2}.
\end{split}
\end{equation}

Next, we treat the terms ${\cal D}_{\omega}$ and ${\cal D}_{\theta}$. Recall that
\begin{equation*}
\begin{split}
{\cal D}_{\theta}= \mu\int_0^T\int_{\Bbb R}\int_{\Bbb T}  A\left((a^2-1)\partial^L_{YY}\theta\right)A\theta dXdYdt.
\end{split}
\end{equation*}
Dividing this term into the zero and nonzero modes gives
\begin{equation*}
\begin{split}
{\cal D}_{\theta}= & \mu\int_0^T\int_{\Bbb R}\int_{\Bbb T}  A\left((a^2-1)\partial^L_{YY}(\theta_0+\theta_{\neq}) \right)A(\theta_0+\theta_{\neq})dXdYdt\\
= & \mu\int_0^T\int_{\Bbb R}\int_{\Bbb T}  A\left((a^2-1)\partial_{YY}\theta_0 \right)A\theta_0dXdYdt  +\mu\int_0^T\int_{\Bbb R}\int_{\Bbb T}  A\left((a^2-1)\partial^L_{YY}\theta_{\neq} \right)A\theta_{\neq}dXdYdt\\
= & {\cal D}_{\theta}^{(0,0)}+{\cal D}_{\theta}^{(\neq,\neq)}.
\end{split}
\end{equation*}
We note that the $(0,\neq)$ and $(\neq,0)$ terms are zero due to the vanishing $X$-average of $\theta_{\neq}$ and $a$ is independent of $X$. For the $(0,0)$ term, we have by integration by parts
\begin{equation*}
\begin{split}
{\cal D}_{\theta}^{(0,0)} = & -2\mu\int_0^T\int_{\Bbb R}\int_{\Bbb T}  A\left(b\partial_{Y}\theta_0 \right)A\theta_0dXdYdt -\mu\int_0^T\int_{\Bbb R}\int_{\Bbb T}  A\left((a^2-1)\partial_{Y}\theta_0 \right)A\partial_Y\theta_0dXdYdt \\
\lesssim & \mu\|b\|_{L^2H^N}\|\partial_Y\theta_0\|_{L^2H^N} \|\theta_0\|_{L^{\infty}H^N} +\mu\|a^2-1\|_{L^{\infty}H^N} \|\partial_Y\theta_0\|^2_{L^2H^N}\\
\lesssim & \mu\delta \varepsilon_2^2\mu^{-1}\lesssim \delta \varepsilon_2^2,
\end{split}
\end{equation*}
by assumption \eqref{assumption}. For the $(\neq,\neq)$ term, we have by integration by parts,
\begin{equation*}
\begin{split}
{\cal D}_{\theta}^{(\neq,\neq)}= &  -2\mu\int_0^T\int_{\Bbb R}\int_{\Bbb T}  A\left(b\partial_{Y}^L\theta_{\neq} \right)A\theta_{\neq}dXdYdt -\mu\int_0^T\int_{\Bbb R}\int_{\Bbb T}  A\left((a^2-1)\partial_{Y}^L\theta_{\neq} \right)A\partial_Y^L\theta_{\neq}dXdYdt\\
\lesssim & \mu\|b\|_{L^{\infty}H^N}\|\nabla_L\theta_{\neq}\|_{L^2H^N} \|A\theta_{\neq}\|_{L^2L^2} + \mu\|a^2-1\|_{L^{\infty}H^N} \|\nabla_L\theta_{\neq}\|_{L^2H^N} \|\nabla_LA\theta_{\neq}\|_{L^2L^2}\\
\lesssim & \mu\delta \|\nabla_L\theta_{\neq}\|^2_{L^2H^N} +\mu\delta\|A\theta_{\neq}\|^2_{L^2L^2}\\
\lesssim & \delta \varepsilon_2^2.
\end{split}
\end{equation*}
This immediately gives
\begin{equation}\label{eq18}
\begin{split}
{\cal D}_{\theta} \lesssim & \delta \varepsilon_2^2.
\end{split}
\end{equation}
Similar estimate leads to
\begin{equation}\label{eq19}
\begin{split}
{\cal D}_{\omega} \lesssim & \delta \varepsilon_1^2.
\end{split}
\end{equation}

For the linear term ${\cal T}_b$, we have
\begin{equation*}
\begin{split}
{\cal T}_b= & (\mu-\nu)\int_0^T\int_{\Bbb R}\int_{\Bbb T}  A(b\partial_Y^L(\theta_0+\theta_{\neq})) A(\theta_0+\theta_{\neq})dXdYdt \\
= & (\mu-\nu)\int_0^T\int_{\Bbb R}\int_{\Bbb T}  A(b\partial_Y^L\theta_0) A\theta_0dXdYdt +(\mu-\nu)\int_0^T\int_{\Bbb R}\int_{\Bbb T}  A(b\partial_Y^L\theta_{\neq}) A\theta_{\neq}dXdYdt\\
= & {\cal T}_b^{(0,0)}+{\cal T}_b^{(\neq,\neq)},
\end{split}
\end{equation*}
where again we have used the fact that the $(0,\neq)$ and $(\neq,0)$ terms are zero due to the vanishing $X$-average of $\theta_{\neq}$ and $a$ is independent of $X$. Similar to the estimates of $D_{\theta}^{(0,0)}$ and $D_{\theta}^{(\neq,\neq)}$,
\begin{equation*}
\begin{split}
{\cal T}_b^{(0,0)}+{\cal T}_b^{(\neq,\neq)}\leq & |\mu-\nu|\delta\varepsilon_2^2\mu^{-1}\\
\lesssim & \delta\varepsilon_2^2,
\end{split}
\end{equation*}
if we assume that $\nu\lesssim \mu$ by a universal constant. In particular, when $\nu=\mu$, this term vanishes.

Below, we treat ${\cal T}_{\omega\theta}$ and ${\cal T}_{\theta\omega}$. By H\"older inequality  and Poincar\'{e} inequality, we have
\begin{equation*}
\begin{split}
{\cal T}_{\omega\theta}=& \int_0^T\int_{\Bbb R}\int_{\Bbb T}  A\omega_{\neq}A\partial_X\theta dXdYdt\leq \|A\omega_{\neq}\|_{L^2L^2}\|\partial_XA\theta\|_{L^2L^2}\\
\leq & \|\nabla_LA\omega_{\neq}\|_{L^2L^2}\|\nabla_LA\theta\|_{L^2L^2}\leq {\varepsilon_1\varepsilon_2}{{\nu^{-1/2}\mu^{-1/2}}},
\end{split}
\end{equation*}
and
\begin{equation*}
\begin{split}
\alpha{\cal T}_{\theta\omega}=& \alpha\int_0^T\int_{\Bbb R}\int_{\Bbb T}  A(\partial_X\Delta_t^{-1}\omega)A\theta_{\neq} dXdYdt \leq \alpha\|A\partial_X\Delta_t^{-1}\omega_{\neq}\|_{L^2L^2} \|A\theta_{\neq}\|_{L^2L^2}\\
\leq & \alpha\|\nabla_LA\theta_{\neq}\|_{L^2L^2} \|A\omega_{\neq}\|_{L^2L^2}\leq \frac{\alpha\varepsilon_1\varepsilon_2} {\mu^{1/2}\nu^{1/6}},
\end{split}
\end{equation*}
where we have used $\|\partial_X\Delta_t^{-1}\|\lesssim 1$. If we assume $\varepsilon_2\ll \sqrt{\nu\mu}\varepsilon_1$ and $0\leq \alpha<\mu^{1/2}\nu^{1/6}\varepsilon_2/\varepsilon_1$, then we have readily
\begin{equation*}
{\cal T}_{\omega\theta}\ll \varepsilon_1^2,\ \ \ \alpha{\cal T}_{\theta\omega}\ll \varepsilon_2^2.
\end{equation*}
From \eqref{eq20} and \eqref{eq21}, under the assumptions $\varepsilon_1\ll\nu^{1/3}$, $\varepsilon_2\ll\sqrt{\nu\mu}\varepsilon_1$, the assumptions \eqref{assumption} indeed hold with `8' replaced by `4'.

Finally, we will prove the third inequality in \eqref{eq8}. We also assume \eqref{assumption} and \eqref{eq2} hold. In particular,
\begin{equation}\label{eq28}
\|u^X_0\|_{L^{\infty}L^2} +\nu^{1/2}\|\partial_Yu^X_0\|_{L^2L^2} \leq 8\varepsilon_1,
\end{equation}
and we need to prove that the constant $8\varepsilon_1$ can be replaced with $4\varepsilon_1$. Note that the $x$-average of the velocity $u^X_0$
satisfy the following equation
$$\partial_t\tilde u^X_0+(\tilde u\cdot\nabla\tilde u^X)_0-\nu\Delta\tilde u^X_0=0,$$
where we have used $\tilde u^Y_0=0$ thanks to the divergence free condition. This is exactly the same equation for the $x$-average of the first velocity component of the Navier-Stokes equation studied in \cite{BVW18}, and the assumption \eqref{eq28} indeed holds with $8\varepsilon_1$ replaced with $4\varepsilon_1$ by the same estimates as in \cite{BVW18}.

\subsection{Proof of Theorem \ref{thm2}}
In the following, we consider the case when $\alpha$ and $\mu$ is priori fixed, not necessarily small. This is also proved by continuity method. We set the following bootstrap assumption
\begin{equation*}
\begin{split}
\alpha\|A\omega(T)\|^2_{L^2}  +\|\nabla_LA\theta\|^2_{L^2} & + \int_0^T\left(\nu\alpha\|\nabla_LA\omega\|_{L^2}^2 +\frac{\alpha}{2}\|\sqrt{-\dot MM}\omega\|^2_{H^N}\right)dt  \\
& + \int_0^T\left(\|\sqrt{-\Delta_L}\sqrt{-\dot MM}\theta\|_{H^N}  +\frac{\mu}{4}\|\Delta_LA\theta\|^2\right)dt \leq 8\varepsilon^2.
\end{split}
\end{equation*}
Similar to the previous estimates, we first consider the time derivative of $\alpha\|A\omega\|^2_{L^2}$, from which one has
\begin{equation}\label{eq12}
\begin{split}
\frac{1}{2}\frac{d}{dt} \alpha\|A\omega\|^2_{L^2}  & + \nu\alpha\|\nabla_LA\omega\|_{L^2}^2 +\alpha\|\sqrt{-\dot MM}\omega\|^2_{H^N} \\
 = & - \alpha\int_{\Bbb R}\int_{\Bbb T}  A(u\cdot\nabla_t\omega)A\omega dXdY + \alpha\int_{\Bbb R}\int_{\Bbb T}  A(b\partial_X\psi)A\omega dXdY\\
&  +  \alpha\nu\int_{\Bbb R}\int_{\Bbb T} A((a^2-1)\partial^L_{YY}\omega)A\omega dXdY +\alpha\int_{\Bbb R}\int_{\Bbb T}  A\partial_X\theta A\omega dXdY.
\end{split}
\end{equation}
For these terms, we have
\begin{equation*}
\begin{split}
\left|\alpha\int_0^T\int_{\Bbb R}\int_{\Bbb T}  A(u\cdot\nabla_t\omega)A\omega dXdYdt\right|\leq \varepsilon_1^3\alpha^{-1/2}\nu^{-1/2} +\varepsilon_1^3\alpha^{-1/2}\nu^{-1/3},
\end{split}
\end{equation*}
and since $b=b(t,Y)$ independent of $X$, by convolution inequality and \eqref{eq3} we have
\begin{equation*}
\begin{split}
\Bigg|\alpha\int_0^T\int_{\Bbb R}\int_{\Bbb T}   A(b\partial_X\psi) & A\omega  dXdYdt\Bigg|=   \left|\alpha\int_0^T\int_{\Bbb R}\int_{\Bbb T}  A(b\partial_X\Delta_{t}^{-1}\omega_{\neq}) A\omega_{\neq}  dXdYdt\right|\\
\leq &\|b\|_{H^{N+2}}\|\sqrt{-\dot MM}\langle D\rangle^N\Delta_L\Delta_t^{-1}\omega_{\neq} \|_{L^2L^2}\|\sqrt{-\dot MM}\langle D\rangle^N\omega_{\neq}\|_{L^2L^2}\\
\leq & \delta \|\sqrt{-\dot MM}\langle D\rangle^N\omega_{\neq}\|^2_{L^2L^2},
\end{split}
\end{equation*}
which can be absorbed by the left of \eqref{eq12} when $\delta$ is small and $\alpha$ is fixed. Similar treatment as in ${\cal D}_{\omega}$ above leads to
\begin{equation*}
\begin{split}
\left|\alpha\nu\int_0^T\int_{\Bbb R}\int_{\Bbb T}  A((a^2-1)\partial^L_{YY}\omega)A\omega dXdYdt \right|\leq \delta\alpha\varepsilon_1^2.
\end{split}
\end{equation*}
Integrating \eqref{eq12} in time, we have from these estimates that
\begin{equation}\label{eq22}
\begin{split}
\alpha\|A\omega(T)\|^2_{L^2}   + \int_0^T \left( \nu\alpha\|\nabla_LA\omega\|^2_{L^2} +\frac{\alpha}{2}\|\sqrt{-\dot MM}\omega\|^2_{H^N}\right)dt& \\
-\alpha\int_0^T\int_{\Bbb R}\int_{\Bbb T}  A\partial_X\theta A\omega dXdYdt -\alpha\|A\omega_{in}\|^2_{H^N} & \lesssim  \varepsilon_1^3\alpha^{-1/2}\nu^{-1/2} +\delta\alpha\varepsilon_1^2.
\end{split}
\end{equation}

Similarly, by taking $A$ to the second equation and taking inner product with $-\Delta_LA\theta$, we obtain
\begin{equation}\label{eq13}
\begin{split}
\frac{d}{dt}\langle A\theta,-\Delta_LA\theta\rangle & +\|\sqrt{-\Delta_L}\sqrt{-\dot MM}\theta\|_{H^N}  +\mu\|\Delta_LA\theta\|^2 \\
= & -\mu\langle A((a^2-1)\partial^L_{YY}\theta),\Delta_LA\theta \rangle - (\mu-\nu)\langle A(b\partial_Y^L\theta),\Delta_LA\theta\rangle\\
&+\alpha \langle A\theta, \Delta_LA\partial_X\Delta_L^{-1}\omega \rangle +\langle A (u\cdot\nabla_t\theta),\Delta_LA\theta\rangle -2\langle \partial_X\partial_Y^LA\theta ,A\theta\rangle,
\end{split}
\end{equation}
where $\langle\cdot,\cdot\rangle=\int dXdY$. By taking Fourier transform and using Young's inequality, we have
\begin{equation*}
\begin{split}
\left|2\langle \partial_X\partial_Y^LA\theta ,A\theta\rangle\right| = & \left|2\int_{\Bbb R}\int_{\Bbb T}\partial_X\partial_Y^LA\theta_{\neq}A\theta_{\neq} dXdY\right|\\
= & 2\left|\sum_{|k|\neq0}\int_{\Bbb R} k(\xi-kt)M^2(k,\xi)\left(1+k^2+\xi^2\right)^N |\widehat\theta_{\neq}|^2d\xi\right|\\
\leq & \sum_{|k|\neq0}\int_{\Bbb R} \left(k^2+(\xi-kt)^2\right)M^2(k,\xi) \left(1+k^2+\xi^2\right)^{N} |\widehat\theta|^2d\xi \\
\leq &\|\Delta_LA\theta\|^2_{L^2},
\end{split}
\end{equation*}
which can be absorbed by the left hand side provided that $\mu\geq 2$. For the convective term, by divergence free condition we have $u_0^Y=0$,  therefore
\begin{equation*}
\begin{split}
\int_0^T\int_{\Bbb R}\int_{\Bbb T}  A (u\cdot\nabla_t\theta)\Delta_LA\theta dXdYdt= &\int_0^T \int_{\Bbb R}\int_{\Bbb T}  A (u^X_0\partial_X\theta)\Delta_LA\theta dXdYdt \\
& +\int_0^T\int_{\Bbb R}\int_{\Bbb T}  A (\nabla^{\bot}_t\Delta_t^{-1}\omega_{\neq}\cdot\nabla_t\theta)\Delta_LA\theta dXdYdt\\
 =& :\tilde{\cal T}_{\theta,0} +\tilde{\cal T}_{\theta,\neq}.
\end{split}
\end{equation*}
For the first term, we have
\begin{equation*}
\begin{split}
\left|\tilde{\cal T}_{\theta,0}\right|= & \int_0^T\int_{\Bbb R}\int_{\Bbb T}  A (u^X_0\partial_X\theta_{\neq})\Delta_LA\theta_{\neq} dXdYdt\\
= & \int_0^T\int_{\Bbb R}\int_{\Bbb T}  \left(\nabla_LA (u^X_0\partial_X\theta_{\neq})-u^X_0\nabla_LA \partial_X\theta_{\neq}\right)\nabla_LA\theta_{\neq} dXdYdt,
\end{split}
\end{equation*}
which, as in ${\cal T}_{\theta,0}$, can be bounded by
\begin{equation*}
\begin{split}
\left|\tilde{\cal T}_{\theta,0}\right|\leq \|\nabla_L\theta\|^2_{L^2H^N} \|\omega\|_{L^\infty H^N} \lesssim \varepsilon_1\varepsilon_2^2\alpha^{-1/2}\mu^{-1}.
\end{split}
\end{equation*}
For the second term, we have
\begin{equation*}
\begin{split}
\left|\tilde{\cal T}_{\theta,\neq}\right|\leq & \|\nabla_t\theta\|_{L^{\infty}H^N} \|\nabla^{\bot}_t\Delta_t^{-1}\omega_{\neq}\|_{L^2H^N} \|\Delta_LA\theta\|_{L^2L^2}\\
\lesssim & \|\nabla_L\theta\|_{L^{\infty}H^N} \|\sqrt{-\dot MM}\omega_{\neq}\|_{L^2H^N} \|\Delta_LA\theta\|_{L^2L^2}\\
\lesssim  & \varepsilon_1\varepsilon_2^2\mu^{-1},
\end{split}
\end{equation*}
where in the second inequality we have used implicitly
$$\frac{\sqrt{k^2+|\xi-kt|^2}}{k^2+|\xi-kt|^2}\leq \sqrt{\frac{-\dot M}{M}}\lesssim \sqrt{-\dot MM},\ \ \ \forall\, |k|\neq0.$$
For the term involving $a,b$, we have
\begin{equation*}
\begin{split}
\left| \mu \int_{\Bbb R}\int_{\Bbb T}A((a^2-1)\partial^L_{YY}\theta)\Delta_LA\theta dXdY\right| \leq & \frac{\mu\delta}{2}\|\Delta_LA\theta\|_{L^2}^2+ \frac{\mu}{2\delta}\|A((a^2-1)\partial^L_{YY} \theta)\|_{L^2}^2\\
\leq & \frac{\mu\delta}{2}\|\Delta_LA\theta\|_{L^2}^2+ \frac{\mu}{2\delta}\|A(a^2-1)\|_{L^2}^2\|A \partial^L_{YY}\theta\|_{L^2}^2\\
\leq & {\mu\delta}\|\Delta_LA\theta\|_{L^2}^2,
\end{split}
\end{equation*}
which can be absorbed when $\delta$ is small. Similarly, we can prove
\begin{equation*}
\begin{split}
\left|(\mu-\nu)\int_{\Bbb R}\int_{\Bbb T} A(b\partial_Y^L\theta)\Delta_LA\theta dXdY\right| \leq &  \frac{|(\mu-\nu)|}{2} \left( \delta\|\Delta_LA\theta\|_{L^2}^2 +\frac{1}{\delta}\|A(b\partial_Y^L\theta)\|_{L^2}^2\right)\\
\leq &  \delta{\left|(\mu-\nu)\right|}\|\Delta_LA\theta\|_{L^2}^2,
\end{split}
\end{equation*}
where we have used the fact that $|\xi-kt|\lesssim (k^2+|\xi-kt|^2)$ in Fourier space for $k\neq0$. When $\delta$ is small and $\nu\leq 2\mu$ (which is obvious since we assume $\mu\geq2$ when $\alpha>0$ is not small), this term can be absorbed by the left hand side. Therefore, we have from \eqref{eq13} after integration in time that
\begin{equation}\label{eq23}
\begin{split}
\|\nabla_LA\theta(T)\|^2_{L^2} &   +\int_0^T\left(\|\sqrt{-\Delta_L}\sqrt{-\dot MM}\theta\|_{H^N}  +\frac{\mu}{4}\|\Delta_LA\theta\|_{L^2}^2\right)dt \\
\lesssim & \|\nabla_LA\theta_{in}\|^2_{L^2}  +\alpha \int_0^T \int_{\Bbb R}\int_{\Bbb T}  A\theta \Delta_LA\partial_X\Delta_L^{-1}\omega dXdYdt +\varepsilon_1\varepsilon_2^2(1+\alpha^{-1/2}) \mu^{-1}.
\end{split}
\end{equation}

Adding together the inequalities \eqref{eq22} and \eqref{eq23}, and noticing that the $\alpha$-terms cancel by integration by parts in $X$-variable, we readily obtain that
\begin{equation}\label{eq14}
\begin{split}
\alpha\|A\omega(T)\|^2_{L^2}  & +\|\nabla_LA\theta\|^2_{L^2}+ \int_0^T\left(\nu\alpha\|\nabla_LA\omega\|_{L^2}^2 +\frac{\alpha}{2}\|\sqrt{-\dot MM}\omega\|^2_{H^N}\right)dt\\
& +\int_0^T\left(\|\sqrt{-\Delta_L}\sqrt{-\dot MM}\theta\|_{H^N}  +\frac{\mu}{4}\|\Delta_LA\theta\|_{L^2}^2\right)dt\\
\lesssim & \alpha\|A\omega_{in}\|^2_{H^N}   +\|\nabla_LA\theta_{in}\|^2_{L^2} + \varepsilon_1^3\alpha^{-1/2}\nu^{-1/2} +\delta\alpha\varepsilon_1^2 +\varepsilon_1\varepsilon_2^2 (1+\alpha^{-1/2})\mu^{-1}.
\end{split}
\end{equation}
When $\alpha$ is fixed, $\mu\geq2$ and $\nu\in (0,1]$, by selecting the initial data $\omega_{in}$ and $\theta_{in}$ appropriately small, sayp
\begin{equation*}
\begin{split}
\alpha\|A\omega_{in}\|^2_{L^2}= \varepsilon_1^2\leq \varepsilon^2 \leq \gamma_1\nu,\ \ \ \|\nabla_LA\theta_{in}\|^2_{L^2}= \varepsilon_2^2\leq \varepsilon^2 \leq \gamma_2\nu,
\end{split}
\end{equation*}
for some $\gamma_1$ and $\gamma_2$ sufficiently small, we get for $\delta$ sufficiently small that
\begin{equation*}
\begin{split}
\alpha\|A\omega(T)\|^2_{L^2}  +\|\nabla_LA\theta\|^2_{L^2} & + \int_0^T\left(\nu\alpha\|\nabla_LA\omega\|_{L^2}^2 +\frac{\alpha}{2}\|\sqrt{-\dot MM}\omega\|^2_{H^N}\right)dt  \\
& + \int_0^T\left(\|\sqrt{-\Delta_L}\sqrt{-\dot MM}\theta\|_{H^N}  +\frac{\mu}{4}\|\Delta_LA\theta\|_{L^2}^2\right)dt \leq 4\varepsilon^2.
\end{split}
\end{equation*}
Finally, \eqref{eq26} is a direct result of \eqref{eq27} and \eqref{eq15}. This concludes the proof by continuity method.


\bigskip
\noindent {\bf Acknowledgments.}
D. Bian is supported by NSFC under the contract 11871005. X. Pu is supported by NSFC under the contract 11871172 and Natural Science Foundation of Guangdong Province of China under 2019A1515012000.
\small

\begin{center}

\end{center}

\begin{thebibliography}{99}
\addcontentsline{toc}{section}{References} {\small

\bibitem{BGM15a} J. Bedrossian, P. Germain, and N. Masmoudi, Dynamics near the subcritical transition of the 3D Couette flow I: Below threshold, {\it Mem. of the AMS}, 266(1294), (2020)v+158.

\bibitem{BGM15b} J. Bedrossian, P. Germain, and N. Masmoudi, Dynamics near the subcritical transition of the 3D Couette flow II: above threshold, arXiv:1506.03721.

\bibitem{BGM17} J. Bedrossian, P. Germain, and N. Masmoudi, On the stability threshold for the 3D Couette flow in Sobolev regularity, \emph{Ann. Math.}, 185, (2017)541-608.

\bibitem{BGM19} J. Bedrossian, P. Germain, and N. Masmoudi, Stability of the Couette flow at high Reynolds number in two dimensions and three dimensions, \emph{Bull. Amer. Math. Soc.}, 56(3), (2019)373-414.

\bibitem{BM13}  J. Bedrossian and N. Masmoudi, Inviscid damping and the asymptotic stability of planar shear flows in the 2D Euler equations, \emph{Publ. Math. l'IHES.}, 122(1), (2013)193-300.

\bibitem{BMV16} J. Bedrossian, N. Masmoudi and V. Vicol, Enhanced dissipation and inviscid damping in the inviscid limit of the Navier-Stokes equations near the 2D Couette flow, \emph{Arch. Ration. Mech. Anal.}, 216(3), (2016)1087-1159.

\bibitem{BVW18} J. Bedrossian, V. Vicol and F. Wang, The Sobolev stability threshold for 2D shear flows near Couette, \emph{J. Nonl. Sci.}, 28, (2018)2051-2075.

\bibitem{BZD2005} R. Bianchini, M.C. Zelati and M. Dolce, Linear inviscid damping for shear flows near Couette in the 2D stably stratified regime, arXiv:2005.09058v1, 2020.

\bibitem{CD81} J.R. Cannon and E. Di Benedetto, The initial problem for the Boussinesq equations with data in $L^p$, \emph{Lecture Notes in Mathematics}, Vol. {771}, Springer: Berlin, 1980.

\bibitem{Chae06} D. Chae, Global regularity for the 2D Boussinesq equations with partial viscosity terms, \emph{Adv. Math.}, 203(2), (2006)497-513.

\bibitem{CN97} D. Chae and H.-S. Nam, Local existence and blow-up criterion for the Boussinesq equations, \emph{Proc. R. Soc. Edinb. Sect. A}, {127}, (1997)935-946.

\bibitem{DM18} Y. Deng and N. Masmoudi, Long time instability of Couette flow in low Gevrey spaces, arXiv:1803.01246v1, 2018.

\bibitem{DWZ2004} W. Deng, J. Wu and P. Zhang, Stability of Couette flow for 2D Boussinesq system with vertical dissipation, arXiv:2004.09292v1, 2020.

\bibitem{DM18} Y. Deng and N. Masmoudi, Long time instability of the Couette flow in low Gevrey spaces, arXiv:1803.01246v1, 2018.


\bibitem{DWZZ18} C.R. Doering, J. Wu, K. Zhao and X. Zheng, Long time behavior of the two-dimensional Boussinesq equations without buoyancy diffusion, \emph{Physica D: Nonlinear Phenomena}, 376, (2018)144-159.


\bibitem{Gold} S. Goldstein, On the stability of superposed streams of fluids of different densities, \emph{Proc. R. Soc. Lond. A}, 132(820), (1931)524-548.

\bibitem{G-G-T16} E. Grenier, Y. Guo, T. T. Nguyen, Spectral instability of general symmetric shear flows in a two-dimensional channel,  \emph{Adv. Math.}, 292, (2016) 52-110.


\bibitem{HL05} T.Y. Hou and C. Li, Global well-posedness of the viscous Boussinesq equations, \emph{Disc. Cont. Dyn. Sys.}, {12} (2005), 1-12.

\bibitem{Jia20} A.D. Ionescu and H. Jia, Inviscid damping near the Couette flow in a channel, \emph{Comm. Math. Phys.}, 374(3), (2020)2015-2096.

\bibitem{Jia20SIAM} H. Jia, Linear inviscid damping near monotone shear flows, \emph{SIAM J. Math. Anal.}, 52(1), (2020)623-652.


\bibitem{LZ11} Z. Lin and C. Zeng, Inviscid dynamical structures near Couette flow, \emph{Arch. Ration. Mech. Anal.}, 200, (2011)1075-1097.

\bibitem{MSZ20} N. Masmoudi, B. Said-Houari and W. Zhao, Stability of Couette flow for 2D Boussinesq system without thermal diffusivity, arXiv:2010.01612v1, 2020.


\bibitem{Reddy} S. Reddy, P. Schmid, J. Baggett and D. Henningson, On stability of streamwise streaks and transition thresholds in plane channel flows, \emph{J. Fluid Mech.}, 365, (1998)269-303.

\bibitem{Rom73} V.A. Romanov, Stability of plane-parallel Couette flow, Funk. Anal. i. Prilozen, 7(1973)62-73.

\bibitem{Synge} J.L. Synge, The stability of heterogeneous liquids, \emph{Trans. Royal Soc. Canada}, 1993.

\bibitem{TW19} L. Tao and J. Wu, The 2d Boussinesq equations with vertical dissipation and linear stability of shear flows, \emph{J. Differential Equations}, 267(3)(2019)1731-1747.

\bibitem{Tay} G.I. Taylor, Effect of variation in density on the stability of superposed streams of fluid, \emph{Proc. Royal Society London. A.}, 132(820):499-523, 1931.

\bibitem{Vil09} C. Villani, Hypocoercivity, \emph{Mem. Amer. Math. Soc.}, 202(950), (2009)iv+141.

\bibitem{WZ18} D. Wei and Z. Zhang, Transition threshold for the 3D Couette flow in Sobolev space, arXiv:1803.01359v1, 2018

\bibitem{WZZ18} D. Wei, Z. Zhang and W. Zhao,  Linear inviscid damping for a class of monotone shear flow in Sobolev spaces, \emph{Comm. Pure Appl. Math.}, 71(4), (2018)617-687.

\bibitem{WZZ20} D. Wei, Z. Zhang and W. Zhao, Linear inviscid damping and enhanced dissipation for the Kolmogrov flow, \emph{Adv. Math.}, 362, (2020)106963.

\bibitem{YL18} J. Yang and Z. Lin, Linear inviscid damping for Couette flow in stratified fluid, \emph{J. Math. Fluid Mech.}, 20(2), (2018)445-472.

\bibitem{Zill17} C. Zillinger, Linear inviscid damping for monotone shear flows, \emph{Trans. Amer. Math. Soc.}, 369(12), (2017)8799-8855.

\bibitem{Zill20} C. Zillinger, On enhanced dissipation for the Boussinesq equations, arXiv.2004.08125v1, 2020.

\bibitem{Zill2011} C. Zillinger, On the Boussinesq equation with non-monotone temperature profiles, arXiv:2011.02316v1, 2020.

}
\end{thebibliography}
\end{document}